\documentclass[12pt]{article}

\usepackage{epsfig}
\usepackage{latexsym}
\usepackage{caption}
\usepackage{amsmath}
\usepackage{amsfonts}
\usepackage{amssymb}
\usepackage{graphicx}
\usepackage{xcolor}

\usepackage{color}

\makeatletter

\newcommand{\Rmnum}[1]{\expandafter\@slowromancap\romannumeral #1@}
\makeatother

\makeatletter
\let\@fnsymbol\@arabic
\makeatother

\begin{document}

\newtheorem{theorem}{Theorem}[section]
\newtheorem{observation}[theorem]{Observation}
\newtheorem{corollary}[theorem]{Corollary}
\newtheorem{algorithm}[theorem]{Algorithm}
\newtheorem{definition}[theorem]{Definition}
\newtheorem{guess}[theorem]{Conjecture}
\newtheorem{claim}[theorem]{Claim}
\newtheorem{problem}[theorem]{Problem}
\newtheorem{question}[theorem]{Question}
\newtheorem{lemma}[theorem]{Lemma}
\newtheorem{proposition}[theorem]{Proposition}
\newtheorem{fact}[theorem]{Fact}

\makeatletter
  \newcommand\figcaption{\def\@captype{figure}\caption}
  \newcommand\tabcaption{\def\@captype{table}\caption}
\makeatother

\newtheorem{acknowledgement}[theorem]{Acknowledgement}

\newtheorem{axiom}[theorem]{Axiom}
\newtheorem{case}[theorem]{Case}
\newtheorem{conclusion}[theorem]{Conclusion}
\newtheorem{condition}[theorem]{Condition}
\newtheorem{conjecture}[theorem]{Conjecture}
\newtheorem{criterion}[theorem]{Criterion}
\newtheorem{example}[theorem]{Example}
\newtheorem{exercise}[theorem]{Exercise}
\newtheorem{notation}[theorem]{Notation}
\newtheorem{solution}[theorem]{Solution}
\newtheorem{summary}[theorem]{Summary}

\newenvironment{proof}{\noindent {\bf
Proof.}}{\rule{3mm}{3mm}\par\medskip}
\newcommand{\remark}{\medskip\par\noindent {\bf Remark.~~}}
\newcommand{\pp}{{\it p.}}
\newcommand{\de}{\em}

\newcommand{\qf}{Q({\cal F},s)}
\newcommand{\qff}{Q({\cal F}',s)}
\newcommand{\qfff}{Q({\cal F}'',s)}
\newcommand{\f}{{\cal F}}
\newcommand{\ff}{{\cal F}'}
\newcommand{\fff}{{\cal F}''}
\newcommand{\fs}{{\cal F},s}

\newcommand{\WReach}{\operatorname{WReach}}
\newcommand{\SReach}{\operatorname{SReach}}
\newcommand{\chitb}{\chi^{\mathrm{tb}}}
\newcommand{\classC}{\mathcal{C}}
\newcommand{\N}{\mathbb{N}}
\newcommand{\G}{\Gamma}
\newcommand{\g}{\gamma}
\newcommand{\wrt}{with respect to }

\newcommand{\ch}{\rm ch}
\newcommand{\ra}{{\rm ran}}
\newcommand{\co}{{\rm col}}
\newcommand{\sco}{{\rm scol}}
\newcommand{\wc}{{\rm wcol}}
\newcommand{\dc}{{\rm dcol}}
\newcommand{\ar}{{\rm arb}}
\newcommand{\ma}{{\rm mad}}
\newcommand{\di}{{\rm dist}}
\newcommand{\dist}{{\rm dist}}
\newcommand{\tw}{{\rm tw}}
\newcommand{\col}{{\rm col}}
\newcommand{\wcol}{{\rm wcol}}
\newcommand{\edp}[2]{#1^{[\natural #2]}}
\newcommand{\epp}[2]{#1^{\natural #2}}
\newcommand{\ind}{{\rm ind}}
\newcommand{\td}{{\rm td}}
\newcommand{\tbw}{{\rm tbw}}
\newcommand{\std}{{\rm std}}
\newcommand{\clos}{{\rm clos}}

\newcommand{\mad}{{\rm mad}}


\title{Low treebandwidth colorings, $p$-centred colorings and generalized
coloring numbers}

\renewcommand{\thefootnote}{\arabic{footnote}}
\author{
	Wenwen Li \footnotemark\footnotemark[1]$~^,$\footnotemark[2]
	\and
	Daqing Yang \footnotemark[1]$~^,$\footnotemark[3]~~~
}

\footnotetext[1]{School of Mathematical Sciences,
	Zhejiang Normal University, Jinhua, Zhejiang 321004, China}
\footnotetext[2]{E-mail: {\small1433828444\texttt{@qq.com}}.}
\footnotetext[3]{Grant number: NSFC, China 12271489.  E-mail: {\small\texttt{dyang@zjnu.edu.cn}}.}

\maketitle

\begin{abstract}
For a graph $G$ and a positive integer $p$, let $\chi^{tb}_p(G)$ be the minimum number of colors in a vertex coloring such that every subgraph induced
by $i$ color classes has treebandwidth at most $i-1$, for $1\le i\le p$. We establish quantitative relationships between these coloring parameters and generalized
coloring numbers. If $\operatorname{wcol}_{2p}(G)\le m$, then $G$ admits a $p$-centred coloring with at most $[(p+1)m]^m$ colors; in particular, $\chi^{tb}_p(G)\le[(p+1)m]^m$.
Conversely, if $\chi^{tb}_{k+1}(G)\le m$, then $\col_k(G) \le 1+4k^4\binom{m}{k+1}^{2}.$
Consequently, a graph class has bounded expansion if and only if its low treebandwidth coloring numbers are uniformly bounded for every fixed parameter.
This gives a characterization of bounded expansion in terms of colorings defined by bounded ancestor distances in rooted forest layouts.

\end{abstract}


\section{Introduction}
\label{sec:introduction}
Generalized coloring numbers measure reachability with respect to a vertex ordering.
They are closely related to centred colorings, low treedepth colorings, and bounded expansion \cite{NOM2008,NOM2008II,Zhu-BGCN}.
Siebertz \cite{SS-2026} gives a recent survey of these parameters and their applications.

Several authors have studied quantitative bounds for generalized coloring numbers. Grohe et al.~\cite{GKRSS2016} studied these bounds and their applications to neighborhood covers. Kreutzer et al.~\cite{KPRS2016} studied generalized coloring numbers on classes of bounded expansion.
For graphs excluding a fixed minor, van den Heuvel et al.~\cite{HMQRS2016} proved linear bounds for strong coloring numbers and polynomial bounds for weak coloring numbers.

There has also been progress on bounds for centred colorings.
Pilipczuk and Siebertz \cite{PS2018} proved that every proper minor-closed class admits centred colorings with a polynomial number of colors in the parameter.
D\c{e}bski et al.~\cite{DFMS2021} improved several  bounds and obtained polynomial bounds for classes excluding a fixed topological minor.
These results give more precise estimates on particular sparse graph classes.

Treedepth is defined by the minimum height of a rooted forest whose closure contains the graph.
A related restriction bounds the distance between adjacent vertices in the forest rather than the height of the entire forest.
This leads to treebandwidth, studied by Jacob, Lochet, and Paul \cite{JLP2025}.
Unlike treedepth, bounded treebandwidth allows rooted forest layouts of arbitrarily large height.

This paper studies colorings in which the union of a small number of color classes has correspondingly small treebandwidth.
For a positive integer $p$, let $\chi^{tb}_p(G)$ denote the minimum number of colors in a coloring of $G$ such that every subgraph induced by $i$ color classes has treebandwidth at most $i-1$, for every $1\le i\le p$.
Precise definitions are given in Section~\ref{sec:preliminaries}.

The first main result gives a centred coloring bound in terms of a weak coloring number at a radius linear in the centredness parameter.
Specifically, if $\wcol_{2p}(G)\le m$, then $\chi^{tb}_p(G)\le ((p+1)m)^m.$
In fact, the coloring constructed is $p$-centred, which is stronger than the required low treebandwidth condition.

The second main result gives a converse estimate.
If $\chi^{tb}_{k+1}(G)\le m$ and $m\geq k+1$, then $\col_k(G)\le1+4k^4\binom{m}{k}^{2}.$
The proof uses an auxiliary graph, as in Zhu \cite{Zhu-BGCN}.
The forest layouts need not have bounded height. Instead, we bound the distance from the minimum-level vertex of a short graph path to each vertex of that path.

Together with the characterization of bounded expansion by generalized coloring numbers \cite{NOM2008,Zhu-BGCN}, these bounds show that
a graph class has bounded expansion if and only if $\chi^{tb}_p$ is bounded on the class for every positive integer $p$.

The paper is organized as follows. Section~\ref{sec:preliminaries} introduces the relevant parameters and records their basic properties.
Section 3 proves the centred coloring bound.
Section 4 establishes the converse estimate for strong coloring numbers.
Section 5 gives the resulting characterization of bounded expansion.

\section{Preliminaries}
\label{sec:preliminaries}

In this paper, every graph we study is finite and simple. Unless stated otherwise, graphs are nonempty.
The length of a path between two vertices is its number of edges. For a positive integer $m$, write $[m]=\{1,\ldots,m\}$.
For a vertex set $X\subseteq V(G)$, let $G[X]$ denote the subgraph induced by $X$.

\subsection{Generalized coloring numbers}

Let $\Pi(G)$ be the set of all linear orderings of
$V(G)$.
For $L\in\Pi(G)$, write $u<_L v$ if $u$ precedes
$v$ in $L$, and let $G_L$ denote the graph $G$
together with the ordering $L$.

\begin{definition}
Let $k$ be a positive integer and let $L\in\Pi(G)$.

A vertex $u$ is \emph{weakly $k$-accessible} from $v$ with respect to $L$ if $u<_L v$ and there exists a $u$--$v$ path $P$ of length at most $k$
such that $u$ is the minimum vertex of $P$ with respect to $L$.

A vertex $u$ is \emph{$k$-accessible} from $v$ with respect to $L$ if $u<_L v$ and there exists
a $u$--$v$ path $P$ of length at most $k$ such that every internal vertex $z$ of $P$ satisfies $v<_L z$.
\end{definition}

Let $Q_k(G_L,v)$ be the set of vertices that are weakly $k$-accessible from $v$, and let
$R_k(G_L,v)$ be the set of vertices that are $k$-accessible from $v$.
Neither set contains $v$. For convenience, let $Q_k[G_L,v]=Q_k(G_L,v)\cup\{v\}$ and  $R_k[G_L,v]=R_k(G_L,v)\cup\{v\}.$

\begin{definition}
The strong \emph{$k$-coloring number} and the \emph{weak $k$-coloring number} of $G$ are defined by
\[
\col_k(G)
=
1+\min_{L\in\Pi(G)}
\max_{v\in V(G)}|R_k(G_L,v)|
\]
and
\[
\wcol_k(G)
=
1+\min_{L\in\Pi(G)}
\max_{v\in V(G)}|Q_k(G_L,v)|,
\]
respectively.
\end{definition}

We also call $\col_k(G)$ the $k$-coloring number. In particular, $\col_1(G)=\wcol_1(G)=\col(G),$
where $\col(G)$ is the ordinary coloring number, equal to the degeneracy of $G$ plus one.
The following bound is given in \cite{KY,Book-N-O}; see also \cite{SS-2026}.

\begin{lemma}
\label{col-wcol}
For eaqch graph $G$ and positive integer $k$,
\[
\wcol_k(G)\leq \bigl(\col_k(G)\bigr)^k.
\]
\end{lemma}

For a finite graph $G$, both sequences stabilize when the radius is at least $|V(G)|-1$.
Their limiting values satisfy
$$\col(G)=\col_1(G)\leq \col_2(G)\leq \cdots\leq \col_\infty(G)=\tw(G)+1,$$
$$\col(G)=\wcol_1(G)\leq \wcol_2(G)\leq \cdots\leq \wcol_\infty(G)=\td(G).$$
See \cite{NOM2006,NOM2008,Book-N-O,SS-2026}. We recall the definitions of treewidth and treedepth below.

\subsection{Rooted forests and treebandwidth}
A \emph{rooted tree} is a tree with a distinguished vertex called its root.
The \emph{level} of a vertex $v$ is the number of vertices on the path from the root to $v$.
Thus, the root has level one. A \emph{rooted forest} is a disjoint union of rooted trees, and its \emph{height} is the maximum level of its vertices.

A rooted forest $F$ defines a partial order on its vertex set:
$x\leq_F y$ if $x$ is an ancestor of $y$.
Every vertex is considered an ancestor of itself.
The \emph{closure} of $F$, denoted by $\clos(F)$, is the graph on $V(F)$ in which two distinct vertices are adjacent if they are comparable in this partial order.

\begin{definition}
The \emph{treedepth} $\td(G)$ of a graph $G$ is the minimum height of a rooted forest $F$ on $V(G)$ such that $G\subseteq\clos(F).$
\end{definition}

For vertices $u$ and $v$ in the same component of a rooted forest $F$, let $\dist_F(u,v)$ denote their  distance in $F$.
For a nonnegative integer $k$, a vertex $v$ is a \emph{$k$-ancestor} of $u$ if $v$ is an ancestor of $u$ and $\dist_F(u,v)=k$.
It is a \emph{$k^{-}$-ancestor} of $u$ if it is an ancestor of $u$ and $\dist_F(u,v)\leq k$.

\begin{definition}
Let $k$ be a nonnegative integer. The \emph{$k$-closure} of a rooted forest $F$, denoted by $k$-$\clos(F)$, is the graph with vertex set $V(F)$ and edge set
\[
\left\{
uv:
u\neq v,\;
\text{$u$ and $v$ are comparable in $F$},\;
\dist_F(u,v)\leq k
\right\}.
\]
\end{definition}
In particular, $0$-$\clos(F)$ is edgeless.

Thus, adjacent vertices in $k$-$\clos(F)$ are comparable in $F$ and have distance at most $k$ in $F$.

\begin{definition}
A \emph{rooted forest layout} of $G$ is a rooted forest $F$ on $V(G)$ such that the endpoints of every edge of $G$ are comparable in the ancestor order of $F$.

The \emph{treebandwidth} of $G$ is
\[
\tbw(G)
=
\min_F
\max_{uv\in E(G)}\di_F(u,v),
\]
where the minimum is taken over all rooted forest layouts of $G$.
\end{definition}

Equivalently, $\tbw(G)$ is the least nonnegative integer $k$ such that
\[
G\subseteq k\text{-}\clos(F)
\]
for some rooted forest $F$ on $V(G)$.
Thus, $\tbw(G)=0$ if and only if $G$ is edgeless.

This definition is equivalent to the rooted tree definition in \cite{JLP2025}.
Indeed, the components of a forest layout can be joined through their roots to form a rooted tree.
This does not change distances within any original component.
Moreover, no edge of $G$ joins distinct components of the forest layout.

\begin{lemma}
\label{lem:subgraph}
If $H$ is a nonempty subgraph of $G$, then
\[
\tbw(H)\leq\tbw(G).
\]
\end{lemma}

\begin{proof}
Let $F$ be a rooted forest layout of $G$ attaining $\tbw(G)$. Deleting edges does not affect the layout.

To delete a vertex $x$ from $F$, attach each child of $x$ to the parent of $x$, if it exists.
If $x$ is a root, make its children roots. This preserves ancestor relations between the remaining vertices and does not increase distances between comparable vertices.

Repeat this operation for all vertices outside $V(H)$. The resulting forest is a layout of $H$ in which the endpoints of every edge are at distance at most $\tbw(G)$. Hence $\tbw(H)\leq\tbw(G)$.
\end{proof}

A \emph{tree decomposition} of $G$ is a pair $(T,\lambda)$, where $T$ is a tree and $\lambda:V(T)\to 2^{V(G)}$ satisfies the following
conditions:
(1) for every $v\in V(G)$, the set $\{x\in V(T):v\in\lambda(x)\}$ is nonempty and induces a connected subtree of $T$;
(2) for every edge $uv\in E(G)$, some bag
$\lambda(x)$ contains both $u$ and $v$.

The width of $(T,\lambda)$ is
\[
\max_{x\in V(T)}|\lambda(x)|-1.
\]
The \emph{treewidth} $\tw(G)$ is the minimum width of a tree decomposition of $G$.

\begin{proposition}
\label{prop:comparison}
For every nonempty graph $G$,
\[
\tw(G)\leq\tbw(G)\leq\td(G)-1.
\]
\end{proposition}

\begin{proof}
Let $F$ be a rooted forest  of height $\td(G)$ such that $G\subseteq\clos(F)$. Any two comparable vertices of $F$ are at distance at most $\td(G)-1$. Thus,
\[
G\subseteq(\td(G)-1)\text{-}\clos(F),
\]
and $\tbw(G)\leq\td(G)-1$.

For the other inequality, let $b=\tbw(G)$. Let  $F$ be a rooted forest on $V(G)$ such that $G\subseteq b$-$\clos(F)$.
For each $x\in V(F)$, take a bag consisting of $x$ and its ancestors at distance at most $b$.

Each edge of $G$ is contained in the bag of its deeper endpoint.
For each vertex $v$, the nodes whose bags contain $v$ are $v$ and its descendants at distance at most $b$.
These nodes induce a connected subtree of $F$. Join the components of $F$ to form a tree. The bags then give a tree decomposition of $G$.

Each bag has at most $b+1$ vertices. Hence $\tw(G)\leq b$.
\end{proof}

Bounded treebandwidth does not imply bounded treedepth.
For example, paths have treebandwidth at most one, but their treedepth is unbounded.
Thus, a graph of bounded treebandwidth need not have a forest layout of bounded height.
The proof in Section~\ref{sec:converse} uses ancestor distances instead.

It is known {\cite{NOM2006Tree-Width, NOM2008, NOM2008II}} that $tw(G)+1\leq twd(G)+1\leq td(G)\leq tw(G)log_{2}n$, where $n$ is the order of $G$.

\subsection{Low treebandwidth colorings}

\begin{definition}
Let $k$ be a positive integer. A \emph{low treebandwidth coloring of order $k$} of a graph $G$ is a vertex coloring $c:V(G)\to[m]$ such that any $i$ color classes induce a graph of treebandwidth at most $i-1$, for $1\leq i\leq\min\{k,m\}$.
That is,
\[
\tbw\bigl(G[c^{-1}(I)]\bigr)\leq |I|-1
\]
for every nonempty set $I\subseteq[m]$ with $|I|\leq k$ and $c^{-1}(I)\neq\varnothing$. Empty induced subgraphs impose no restriction.

The minimum number of colors in such a coloring is denoted by $\chi^{tb}_k(G)$.

\end{definition}
Equivalently, every nonempty subgraph $H$ of $G$ receives at least $\min\{\tbw(H)+1,k+1\}$ colors.

\begin{observation}
\label{obs:basic-colouring}
For every graph $G$ and positive integer $k$,
\[
\chi^{td}_k(G)\leq\chi^{tb}_{k+1}(G).
\]
Moreover, $\chi^{tb}_1(G)$ is the ordinary chromatic number of $G$, and $\chi^{tb}_2(G)$ is its acyclic chromatic number.
\end{observation}

\begin{proof}
The first assertion follows from the definition.

A graph has treebandwidth zero if and only if it is edgeless.
Thus, a low treebandwidth coloring of order one is a proper coloring.

A graph has treebandwidth at most one if and only if it is a forest.
Indeed, every forest has such a layout, and every graph of treebandwidth at most one is a subgraph of its layout forest.
Hence a low treebandwidth coloring of order two is a proper coloring in which any two color classes induce a forest. This is an acyclic coloring.
\end{proof}

\begin{definition}
A \emph{centred coloring} of a graph $G$ is a
vertex coloring such that every nonempty connected subgraph has a color appearing exactly once.
\end{definition}

\begin{definition}
Let $p$ be a positive integer.
A \emph{$p$-centred coloring} of a graph $G$
is a vertex coloring such that every nonempty
connected subgraph either has a color appearing
exactly once or receives at least $p$ colors.
\end{definition}


In \cite{NOM2006} proved the following lemma.
\begin{lemma}
\label{lemma:p-centerd-td}
Let $G$ be a graph, and $c$ is a $p$-centred coloring of $G$, then any subgraph $H$ of $G$ with $td(H)\leq p$ gets at least $\td(H)$ colors.
\end{lemma}

\begin{lemma}\label{lem:centered-low-tbw}
Let $i$ be a positive integer. If a graph $H$ admits a centred coloring using at most $i$ colors, then $\operatorname{td}(H)\leq i$.
Consequently, for every positive integer $k$, every $k$-centred coloring of a graph $G$ is a low
treebandwidth coloring of order $k$.
\end{lemma}

\begin{proof}
The first assertion is immediate if $H$ is empty. Suppose that $H$ is nonempty, and let $c$ be a centred coloring of $H$ using at most $i$ colors.
Set $p=|V(H)|+1$. Since a centred coloring is $p$-centred for every positive integer $p$, and
\[
\operatorname{td}(H)\leq |V(H)|<p,
\]
Lemma~\ref{lemma:p-centerd-td} implies that $H$ receives at least $\operatorname{td}(H)$ colors under $c$.
Hence
\[
\operatorname{td}(H)\leq |c(V(H))|\leq i.
\]

Now let $c:V(G)\to[m]$ be a $k$-centred coloring of $G$, and let $I\subseteq[m]$ satisfy $1\leq |I|\leq \min\{m, k\}$.
Write $F=G[c^{-1}(I)]$, and suppose that $F$ is nonempty. Every nonempty connected subgraph of $F$ receives at most $|I|\leq k$ colors under $c$.
By the $(k+1)$-centred property, it therefore contains a color that appears on exactly one vertex.
Thus the restriction of $c$ to $V(F)$ is a centred coloring of $F$ using at most $|I|$ colors.
Applying the first assertion gives
\[
\operatorname{td}(F)\leq |I|.
\]
Consequently,
\[
\operatorname{tbw}(F)
\leq \operatorname{td}(F)-1
\leq |I|-1.
\]
Since empty induced subgraphs impose no restriction, $c$ is a low treebandwidth coloring of order $k$.
\end{proof}

\section{Bounding $p$-centred colorings by using weak coloring numbers}
\label{sec:upper}
Zhu (\cite{Zhu-BGCN}, Theorem 2.6) proved that, for every integer $p\geq2$, if the weak $2^{p-2}$-coloring of a graph $G$ is less than or equal to $m$, i.e., wcol$_{2^{p-2}}$($G$)$\leq m$, then the graph $G$ has a $p$-centred coloring using at most $m$ colors. Motivated by this result, we now bound $p$-centred coloring in terms of
$\wcol_{2p}(G)$.
We compare Theorem~\ref{thm:centred} with the result above. To bound the number of colors in a $p$-centred coloring, we use a weak reachability radius of $2p$ instead of $2^{p-2}$.

\begin{theorem}
\label{thm:centred}
If $G$ is a graph with $wcol_{2p}(G)\leq m$, then $G$ has a $p$-centred coloring using at most $((1+p)m)^{m}$.
\end{theorem}

\begin{proof}
Let $G=(V,~E)$ be a graph and $L=v_{1}v_{2}\cdots v_{n}$ is an ordering of $V$ with max$_{v\in V(G)}|Q_{2p}(G_{L},v)|\leq m-1$.
We greedily colored the vertices of the graph $G$ with respect to the order $L$, thus the color assigned to $v$ is different from the color assigned to the vertices in $Q_{2p}(G_{L},~v)$. Since $|Q_{2p}(G_{L},~v)|\leq m-1$, $m$ colors are sufficient. Assume that this coloring process causes each vertex $v$ to be colored with color $f(v)$, and we call $f(v)$ the native color of vertex $v$.

Let $r_{uv}$ be the length of the shortest path of all the paths in graph $G$ witnessing that $u$ is $p$-weakly accessible from $v$ with respect to $L$. In the above coloring process, when a vertex $x$ is colored by $f(x)$, suppose that,
 $$ \textrm{Y} [x]:=\{y_{i}: 1\leq i\leq t, y_{i}\in \textrm{Q}_{p}[G_{L}, v],y_{1}<_{L}y_{2}<_{L}\cdots <_{L}y_{t}=x\}.$$
For each $y_{i}\in \textrm{Y}[x]$, suppose that the shortest length of all the path that witnesses $y_{i}\in \textrm{Q}_{p}[G_{L}, x]\cap \textrm{Y}[x]$ is $r_{i}$. By the definition of $\textrm{Y}[x]$, $r_{i}$ is well defined for $y_{i}$, and $0\leq s_{i}\leq p$.

Now for each color $\alpha$ that has been used in the above coloring process for $f$, we break $\alpha$ into $q+1$ shares, which is indicated by $(\alpha, 0)$, $(\alpha, 1)$, $\cdots$, $(\alpha, q)$.
Then we assign a label $\lambda(x)$ to $x$, where $\lambda(x):=((f(y_{1}), r_{1}), (f(y_{2}), r_{2}), \cdots, (f(y_{t}), r_{t}))$.
Note that we always have $(f(y_{t}), r_{t}))=(f(x), 0)$.
We again color the vertices of $G$, one by one according to the linear ordering $L$ of $V(G)$. Suppose we are about to color vertex $u$, for each to be colored vertex $u$, we shall preassign $u$ a color $\lambda(u)$.
We claim that such a coloring is a $p$-centred coloring of the graph $G$. For any connected subgraph $H$ of a graph $G$, let $v$ be the vertex in $H$ that is a minimum with respect to $L$.
If $\lambda(v)$ appears only once in $H$, the proof is complete. Suppose $\lambda(v)$ appears more than once in $H$. In the following we prove that $H$ uses at least $p$ colors.

Let $u\neq v$ is a vertex in $H$ and $\lambda(u)=\lambda(v)$. Let $P_{0}=(v=v_{0}, v_{1},\cdots, v_{q}=u)$ is the shortest path that connecting $u$ and $v$ in $G$, and where $v_{i}\in V(H)$($0\leq i\leq q$). Then we have $q > 2p$, for otherwise $v$ is weakly $2q$-accessible from $u$, i.e., $v\in Q_{2q}(G_{L}, u)$, thus $f(u)\neq f(v)$,
and since  $(f(u), 0)\in \lambda(u)$, $(f(v), 0)\in \lambda(v)$, this contradicts the assumption of $\lambda(u)=\lambda(v)$.

Fact: There is no edge $e=v_{i}v_{j}$ in $G$, where $v_{i}$, $v_{j}\in V(P_{uv})$ and $|i-j|\geq2$; for otherwise, we replace $P_{0}$ by $P'_{0}=(v=v_{0}v_{1}\cdots v_{i}v_{j}v_{j+1}\cdots v_{q}=u)$ and $||P'_{0}| < ||P_{0}||$, a contradiction.

Let $v_{t+1}\in V(P_{0})$ be the vertex with the smallest index in the linear order $L$ such that such that $r_{v_{0}v_{t+1}} > p$ holds.
By the definition of $P_{0}$ this vertex must exist and $p\leq t <q$.
This means $v_{t}\in V(P_{0})$ in $G$ is weakly $p$-reachable from $v_{0}$, and the shortest length of all paths that witnesses $v_{0}\in Q_{p}(G_{L}, v_{t})$ is $p$, i.e., $r_{v_{0}v_{t}}=p$.
Suppose the contrary that $r_{v_{0}}$ then $r_{v_{0}v_{t+1}}\leq r_{v_{0}v_{t}}+||v_{t}v_{t+1}||\leq ||P_{0}||$, a contradiction. Let $P_{1}=\{v=v_{0}v_{1}\cdots v_{t}\}$, $p\leq t< q$. In the following we show that at least $p$ colors appear on $P_{1}$ is the final coloring.

Note that for any vertex $v_{i}\in V(P_{1})$, $(f(v_{0}),1)\in \lambda(v_{1})$, $(f(v_{0}), r_{v_{0}v_{i}})\in \lambda (v_{i})$ and $(f(v_{0}), p)\in \lambda(v_{t})$. Firstly, we show that for any positive integer $k$, $1\leq k\leq q$, there exists a vertex $v_{i_{k}}\in V(P_{1})$ such that $r_{v_{0}v_{i_{k}}}=k$, i.e.,  $(f(v_{0}), k)\in \lambda(v_{i_{k}})$.

Conversely, suppose there exists an integer $k\leq p$ is the minimum integer such that there exists $k-1$ vertices $\{ v_{i_{1}}=v_{1}, v_{i_{2}}, \cdots, v_{i_{k-1}}\}\subseteq V(P_{1})$ and for each $1\leq j\leq k-1$, $(f(v_{0}), j)\in \lambda(v_{i_{j}})$, and there exists no vertex $v_{i_{k}}\in V(P_{1})$ such that $(f(v_{0}), k)\in \lambda(v_{i_{k}})$.

Let $v_{s}\in \{v_{i_{1}}, v_{i_{2}}, \cdots, v_{i_{k-1}}\}$ is the vertex with maximum index in $\{ v_{i_{1}}, v_{i_{2}}, \cdots, v_{i_{k-1}} \}$ in the linear order $L$, then $(f(v_{0}),r_{v_{0}v_{s}})\in \lambda(v_{s})$, $1\leq r_{v_{0}v_{s}}\leq k-1$.
Which means that the shortest path $P_{v_{0}v_{s}}$ that witnesses $v_{0}\in Q_{p}(G_{L}, v_{s})$ in $G$ has length at most $k-1$.
By the assumptions, there is $(f(v_{0}), r_{v_{0}v_{s+1}})\in \lambda(v_{s+1})$ in the color $\lambda(v_{s+1})$ that is color by $v_{s+1}$,
where $r_{v_{0}v_{s+1}}\geq k+1$, i.e., the shortest path $P_{v_{0}v_{s+1}}$ that witnesses $v_{0}\in Q_{p}(G_{L},v_{s+1})$ in $G$ has length at least $k+1$.
Let $P'_{v_{0}v_{s+1}}=P_{v_{s}v_{s+1}}+v_{s}v_{s+1}$, then $P'_{v_{0}v_{t+1}}$ is also path in $G$ that witnesses $v_{0}\in Q_{p}(G_{L}, v_{s+1})$ and $||P'_{v_{0}v_{s+1}}||=||P_{v_{0}v_{s}}||+1\leq k$, this contradicts $r_{v_{0}v_{s+1}}\geq k+1$.

Let $V^{*}=\{v_{i_{1}}=v_{1}, v_{i_{2}},\cdots, v_{i_{p}}\}\subseteq V(P_{1})$ and $(f(v_{i_{j}}), j) \in \lambda(v_{i_{j}})$ for each integer $j$ ($1\leq j\leq q$).
If there exists $v_{i_{j_{1}}}$, $v_{i_{j_{2}}}\in V^{*}$, $\lambda(v_{i_{j}})=\lambda(v_{i_{j'}})$ and $v_{i_{j}}\neq v_{i_{j'}}$,
then there exists a vertex $w\in V(G)$, such that $(f(w), r_{wv_{i_{j'}}})=(f(v_{0}), r_{v_{0}}v_{i_{j}})\in\lambda(v_{i_{j'}})$,
where $r_{wv_{i_{j'}}}=r_{v_{0}}v_{i_{j}}\leq p$,
thus the walk concatenated by the witnessing paths of $v_{0}, w\in Q_{p}(G_{L}, v_{i_{j'}})$ will contain a path that witnesses $v_{0}\in Q_{2p}(G_{L}, w)$ (if $v_{0} <_{L} w$) or $w\in Q_{2p}(G_{L}, v_{0})$ (if $w <_{L} v_{0}$), in either case, $f(w)\neq f(v_{0})$. This proves that any two points in $V^{*}$ are assigned different colors, which completes the proof of the theorem.
\end{proof}

\begin{corollary}
\label{cor:chi-upper}
For every graph $G$ and positive integer $p$,
\[
\chi^{tb}_p(G)
\leq
\left[
(p+1)\wcol_{2p}(G)
\right]^{\wcol_{2p}(G)}.
\]
\end{corollary}

\begin{proof}
Apply Theorem~\ref{thm:centred} and
Lemma~\ref{lemma:p-centerd-td}.
\end{proof}

Theorem~\ref{thm:centred} applies to every graph. For particular graph classes, direct constructions can give smaller bounds for centred colorings.
Examples include proper minor-closed classes \cite{PS2018} and classes excluding a fixed topological minor \cite{DFMS2021}.
These results use the alternative convention for centredness described in Section~\ref{sec:preliminaries}.

\section{Bounding generalized coloring numbers by using low treebandwidth coloings}
\label{sec:converse}
We next bound $\col_k(G)$ in terms of a low treebandwidth coloring number.
The following property of paths in rooted forest layouts replaces the bounded-height argument used for low treedepth colorings.

\begin{theorem}
\label{thm:converse}
Let $k$ be a positive integer, $G$ is a graph with $\chi^{tb}_{k+1}(G)\leq m$, then $col_{k}(G)\leq 1+[2k^2\binom{m}{k+1}]^{2}$.
\end{theorem}

\begin{proof}
The proof of this theorem is similar to the proof of Theorem 2.10 in \cite{Zhu-BGCN}. Let $c$ be an $m$-coloring of $V(G)$ such that for any subgraph of $G$ of treebandwidth  $i$ gets at least $i+1$ colors. Let $[m]=\{1,2,\cdots, m\}$ be the color set, and let $\binom{m}{k+1}$, $c^{-1}(I)$ induces a subgraph $G_{I}$ of treebandwidth at most $k$. Let $F_{I}$ be the rooted forest of height $k+1$ such that $G_{I}$ is a subgraph of $k$-clos($F_{I}$).

Let $\overrightarrow{G_{I}}$ be the directed graph obtained from $k$-clos($F_{I}$) by orient its edges so that $xy$ is oriented from $x$ to $y$ if $y$ is a $k^{-}$-ancestor of $x$. Thus $\overrightarrow{G_{I}}$ has maximum out-degree at most $k$. Furthermore, let $\overrightarrow{G^{+}_{I}}$ be the directed graph obtained from the digraph $\overrightarrow{G_{I}}$ by adding an arc($x,y$) if $x$ and $y$ are not adjacent and there exists in $\overrightarrow{G_{I}}$ a directed path of length at most $k$ from $x$ to $y$.
Hence the maximum out-degree of $\overrightarrow{G^{+}_{I}}$ is at most $k^{2}$. Let
$$H=\bigcup\limits_{I\in\binom{[m]}{k+1}}\overrightarrow{G^{+}_{I}}.$$

The union of the orientations of all the $\overrightarrow{G^{+}_{I}}$'s induces an orientation of $H$ in which each vertex has out-degree at most $\binom{m}{k+1}k^{2}$ (if an edge $e$ of $H$ appears in more than one of the $k$-clos($F_{I}$)'s, the we arbitrarily choose one of the $k$-clos($F_{I}$) to determine the orientation of $e$). This implies that any subgraph $H'$ of $H$ has average degree at most $2\binom{m}{k+1}k^{2}$. Therefore, $H$ has coloring number at most $2\binom{m}{k+1}k^{2}$+1. Let $M=2\binom{m}{k+1}k^{2}$ and let $L$ be a linear ordering of the vertices of $H$ such that for each vertex $v$, $|R_{k}(G_{L}, v)|\leq M$.

Now we prove that for each vertex $v$, $|R_{k}(G_{L}, v)|\leq M^{2}$. Let $v$ be a fixed vertex of $G$. For each $u\in R_{k}(G_{L}, v)$, choose a $u$-$v$-path $P_{uv}$ of length at most $k$ such that all interior vertices $z$ of $P_{uv}$ satisfy $v<_{L}z$. By definition of
$R_{k}(G_{L}, v)$, such a path exists. As $|V(P_{uv})|\leq k+1$, there is an $I\in\binom{m}{k+1}$ such that $c(P_{uv})\subseteq I$. For each $u\in R_{k}(G_{L},v)$, fix $I_{u}\in \binom{m}{k+1}$ such that $c(P_{uv})\subseteq I_{u}$. Let $z_{u}\in P_{uv}$ be the vertex of $P_{uv}$ which has minimum level in $F_{I_{u}}$. If $z_{u}=v$, then $v$ is a  $k^{-}$-ancestor of $u$ in $F_{I_{u}}$ and thus arc $(u,v)\in \overrightarrow{G^{+}_{I}}$. So we have
$$|\{u\in R_{k}(G_{L},v):z_{u}=v\}|\leq |R_{1}(H_{L},v)|\leq M.$$

Similarly, if $z_{u}=u$, we have $\{u\in R_{k}(G_{L}, v):~z_{u}=u\}\subseteq R_{1}(H_{L}, v). $ Assume $z_{u}\neq v,u$. Then $z_{u}$ is a $k^{-}$-ancestor of $v$ and $u$ in $F_{I_{u}}$. As $v$ has at most $k(k-1)$ adjacent ancestors in $H$ and $u\in R_{1}(H_{L},z_{u})$, we conclude that
$$|\{u\in R_{p}(G_{L},v):z_{u}\neq v\}|\leq \binom{m}{k+1}k(k-1)M\leq \frac{M^{2}}{2}.$$

Hence,
$$|R_{k}(G_{L},v)|\leq M+\frac{M^{2}}{2}<M^{2}.$$

\end{proof}

\section{Low treebandwidth colorings and bounded expansion}
\label{sec:expansion}

A \emph{depth-$r$ shallow minor} of a graph $G$ is obtained from a subgraph of $G$ by contracting pairwise vertex-disjoint connected subgraphs, each of radius at most $r$, and then deleting loops and parallel edges.

A class $\mathcal{C}$ of graphs has \emph{bounded expansion} if, for every
nonnegative integer $r$, there is a constant $d_r$ such that every nonempty depth-$r$ shallow minor $H$ of each graph in
$\mathcal{C}$ satisfies
\[
\frac{|E(H)|}{|V(H)|}\leq d_r.
\]

A class $\mathcal{C}$ has \emph{low treebandwidth colorings} if $\chi_k$ is bounded on $\mathcal{C}$ for every positive integer $k$.

Bounded expansion can be characterized by generalized coloring numbers and by centred colorings \cite{NOM2008,Zhu-BGCN}.
See also \cite{KPRS2016,SS-2026}. We add low treebandwidth colorings to these characterizations.

\begin{theorem}
\label{thm:equivalence}
For a class $\mathcal{C}$ of nonempty finite graphs, the following statements are equivalent.
\begin{enumerate}
\item
$\mathcal{C}$ has bounded expansion.

\item
For every positive integer $k$, there is a
constant $a_k$ such that
\[
\col_k(G)\leq a_k
\qquad\text{for every }G\in\mathcal{C}.
\]

\item
For every positive integer $k$, there is a
constant $b_k$ such that every $G\in\mathcal{C}$
has a $k$-centred coloring using at most
$b_k$ colors.

\item
For every positive integer $k$, there is a
constant $c_k$ such that
\[
\chi^{tb}_k(G)\leq c_k
\qquad\text{for every }G\in\mathcal{C}.
\]
\end{enumerate}
\end{theorem}

\begin{proof}
The equivalence of (1) and (2) is established in \cite{Y-2009,Zhu-BGCN}.
We prove
\[
(2)\Longrightarrow(3)\Longrightarrow(4)
\Longrightarrow(2).
\]

Suppose that (2) holds. Fix a positive integer $p$. Choose a positive integer $a_{2p}$ such that
\[
\col_{2p}(G)\leq a_{2p}
\qquad\text{for every }G\in\mathcal{C}.
\]
By Lemma~\ref{col-wcol},
\[
\wcol_{2p}(G)
\leq
(a_{2p})^{2p}
\qquad\text{for every }G\in\mathcal{C}.
\]
Theorem~\ref{thm:centred} gives a $p$-centred coloring of each such graph using at most
\[
[(p+1)(a_{2p})^{2p}]^{(a_{2p})^{2p}}
\]
colors. This coloring is also $p$-centred.
Hence (3) holds.

Suppose that (3) holds. Fix a positive integer $k$. Every $G\in\mathcal{C}$ has a $k$-centred coloring using at most $b_k$ colors. By Lemma~\ref{lem:centered-low-tbw}, this is a low treebandwidth coloring of order $k$. Thus,
\[
\chi_k(G)\leq b_{k}
\qquad\text{for every }G\in\mathcal{C},
\]
which proves (4).

Finally, suppose that (4) holds. Fix a positive integer $k$, and choose an integer
\[
m\geq\max\{c_{k+1},k+1\}.
\]
Then $\chi^{tb}_{k+1}(G)\leq m$ for every $G\in\mathcal{C}$. By Theorem~\ref{thm:converse},
\[
\col_k(G)
\leq
1+\left[2k^2\binom{m}{k+1}\right]^2
\qquad\text{for every }G\in\mathcal{C}.
\]
The right-hand side is independent of $G$, so (2) holds.
\end{proof}

\noindent\textbf{Concluding Remarks:}
Low treebandwidth colorings characterize bounded expansion, just as low treedepth colorings do.
The quantitative bounds depend on two constructions.
The first records weak-accessibility distances in a color vector.
The second adds edges between vertices and nearby ancestors in forest layouts.
Neither construction requires a bound on the height of a treebandwidth layout.

The dependence on the parameters remains to be determined more precisely.
In particular, it would be useful to improve the bound for $\chi^{tb}_p(G)$ in terms of $\wcol_{2p}(G)$ and the bound for $\col_k(G)$ in terms of $\chi^{tb}_{k+1}(G)$.

Another question is how low treebandwidth and low treedepth coloring numbers compare on particular sparse graph classes.
The centred coloring bounds in \cite{PS2018,DFMS2021} also give bounds for low treebandwidth colorings.
It is natural to ask when the weaker treebandwidth requirement permits substantially fewer colors.

\bibliographystyle{plain}
\bibliography{Ref1}

\end{document}